\documentclass[10pt,reqno]{amsart}

\usepackage{bcol}
\usepackage{hyperref}
\newcommand{\triangleq}{\equiv}
\usepackage{algorithm,algorithmicx}
\usepackage{algpseudocode}
\usepackage{rotating}

\usepackage{graphicx,accents}

\newtheorem{definition}{Definition}

\newtheorem{proposition}{Proposition}

\makeatletter
\DeclareRobustCommand{\cev}[1]{%
	\mathpalette\do@cev{#1}%
}
\newcommand{\do@cev}[2]{%
	\fix@cev{#1}{+}%
	\reflectbox{$\m@th#1\vec{\reflectbox{$\fix@cev{#1}{-}\m@th#1#2\fix@cev{#1}{+}$}}$}%
	\fix@cev{#1}{-}%
}
\newcommand{\fix@cev}[2]{%
	\ifx#1\displaystyle
	\mkern#23mu
	\else
	\ifx#1\textstyle
	\mkern#23mu
	\else
	\ifx#1\scriptstyle
	\mkern#22mu
	\else
	\mkern#22mu
	\fi
	\fi
	\fi
}

\newcommand\ubar[1]{%
	\underaccent{\bar}{#1}}

\title{A Scalable Semidefinite Relaxation Approach to Grid Scheduling}

\author{Ramtin Madani}
\author{Alper Atamt\"{u}rk} 
\author{Ali Davoudi}

\ignore{
\affil[a]{Department of Electrical Engineering, The University of Texas at Arlington, TX 76019}
\affil[b]{Department of Industrial Engineering and Operations Research, University of California, Berkeley, CA 94720}

\leadauthor{Madani}

\significancestatement{
This paper is concerned with the problem of scheduling supply of electricity through available resources in a power grid, to meet the demand within a time horizon (e.g. the day-ahead). Current practice for solving this short-term planning problem considers an approximate model for describing the flow of electricity across the grid, instead of an accurate physical model, whose integration is estimated to save billions of dollars per year, by enhancing the efficiency and reliability of operation \cite{cain2012history}. The primary obstacle towards the realization of this direction is the computational complexity incurred by the accurate model. In this paper, a computational method is developed that is compatible with the accurate model for power grids and offers considerable improvements, in terms of scalability.
}

\authorcontributions{Please provide details of author contributions here.}
\authordeclaration{Please declare any conflict of interest here.}
\equalauthors{\textsuperscript{1}A.O.(Author One) and A.T. (Author Two) contributed equally to this work (remove if not applicable).}
\correspondingauthor{\textsuperscript{2}To whom correspondence should be addressed. E-mail: author.two\@email.com}


\dates{This manuscript was compiled on \today}
\doi{\url{www.pnas.org/cgi/doi/10.1073/pnas.XXXXXXXXXX}}

\begin{document}

\verticaladjustment{-2pt}

\maketitle
\thispagestyle{firststyle}
\ifthenelse{\boolean{shortarticle}}{\ifthenelse{\boolean{singlecolumn}}{\abscontentformatted}{\abscontent}}{}

}

\begin{document}

\begin{abstract}
	Determination of the most economic strategies for supply and transmission of electricity is a daunting computational challenge. Due to theoretical barriers, so-called NP-hardness, the amount of effort to optimize the schedule of generating units and route of power, can grow exponentially with the number of decision variables. Practical approaches to this problem involve legacy approximations and ad-hoc heuristics that may undermine the efficiency and reliability of power system operations, that are ever growing in scale and complexity. Therefore, the development of powerful optimization methods for detailed power system scheduling is critical to the realization of smart grids and has received significant attention recently. In this paper, we propose for the first time a computational method, which is capable of solving large-scale power system scheduling problems with thousands of generating units, while accurately incorporating the nonlinear equations that govern the flow of electricity on the grid. The utilization of this accurate nonlinear model, as opposed to its linear approximations, results in a more efficient and transparent market design, as well as improvements in the reliability of power system operations. We design a polynomial-time solvable third-order semidefinite programming (TSDP) relaxation, with the aim of finding a near globally optimal solution for the original NP-hard problem. The proposed method is demonstrated on the largest available benchmark instances from real-world European grid data, for which provably optimal or near-optimal solutions are obtained. 
\end{abstract}

	\maketitle
	
	\begin{center} July 2017 \end{center}
	
	\BCOLReport{17.03}

	The development of systematic algorithms for optimal allocation and scheduling of resources can be traced back to Nobel Laureate Wassily Leontief's seminal work on input-output analysis \cite{leontief1963multiregional}, and the development of simplex algorithm by George Dantzig \cite{dantzig2016linear}, which was driven by the United States military planning problems during World War II. Ever since, various application domains have relied heavily on optimization theory as a tool for design, planning, and decision-making. 


An economically significant application is the operation of power grids, with the aim of efficient and secure supply of electricity.
Grid operations are currently planned with legacy frameworks that are far from producing near-optimal solutions at the scale and detail required by the next-generation grids.
In the past decade, various potential approaches have been identified for enhancements in grid operation through more accurate models for the flow of power, control of network topology, and taking uncertainties of demand and renewable generation into consideration.
The realization of the aforementioned directions is expected to offer considerable improvements in the efficiency and reliability of the power grids \cite{cain2012history}. However, due to the ever-growing size and scope of grids, the scalability of algorithms for solving detailed and accurate models remains as the primary bottleneck \cite{national2016analytic,clack2017evaluation,giannakis2013monitoring}.



\begin{figure*}
	\centering
	\includegraphics[width=0.9\linewidth]{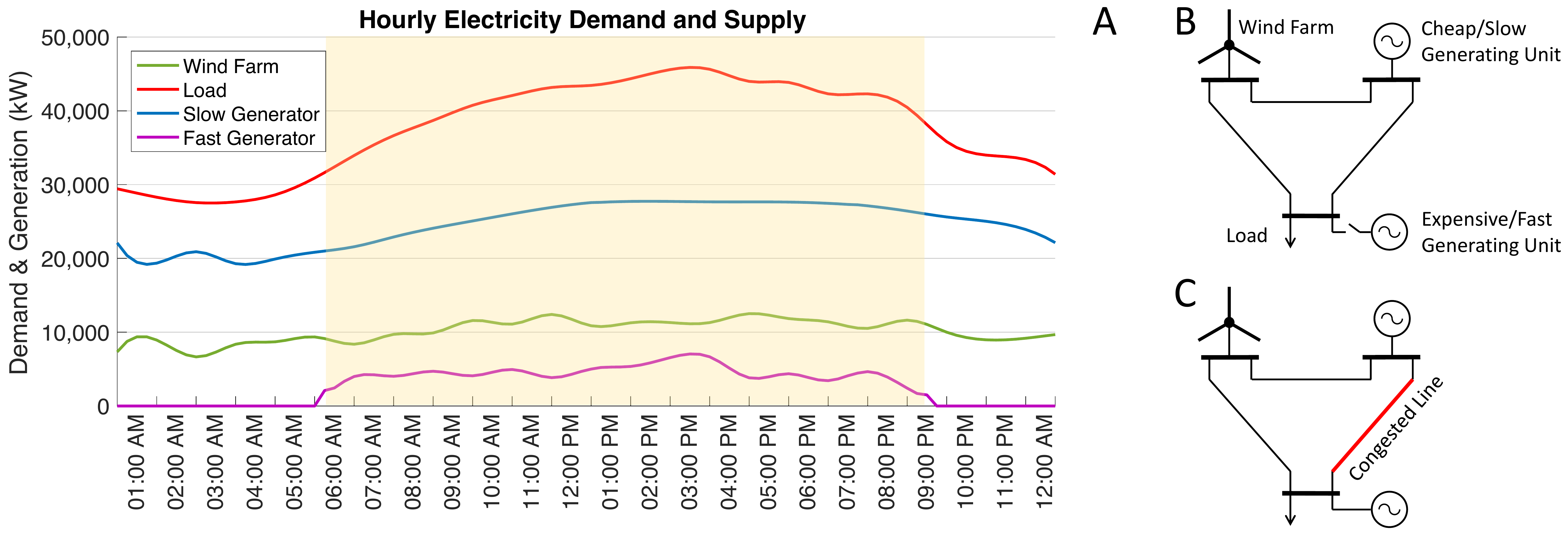}
	\caption{Day-ahead scheduling of a notional power system, with three vertices and two generating units. (A) The optimal operational strategy based on the available forecast of the demand and wind generation. The shaded period represents peak hours. (B) Off-peak configuration of the network, during which the expensive generator is not committed. (C) Peak configuration, during which the expensive generator contributes to accommodate transmission limits.}
	\label{fig:UC_OPF}
\end{figure*}

Building an optimal day-ahead plan for the operation of a nationwide grid is a daunting challenge, in part, due to the presence of thousands of generating units, whose on/off status need to be determined. 
Algorithms for finding the most economic plan with binary on/off decisions give rise to massive search trees.  Another challenge is posed by the nonlinearity of the physical laws describing the flow of electricity.

This paper presents the first computational method that is capable of solving day-ahead power system scheduling problems of realistic size, that are built upon accurate physical models. 
More specifically, the proposed method accomplishes the integrated optimization of two fundamental problems, faced by system operators and utilities on a daily basis: 
\begin{enumerate}
	\item {\bf Unit Commitment (UC):} The problem of scheduling generating units throughout a planning horizon, based on demand forecasts and technological constraints.
	\item {\bf Optimal Power Flow (OPF):} The problem of determining an operating point for the network that delivers power from suppliers to consumers as economically as possible, subject to physical constraints.
\end{enumerate}
Advanced algorithms for UC and OPF can contribute to the efficiency and transparency of power markets by improving operational decisions and pricing mechanisms \cite{lipka2017running,7778217}.
Figure \ref{fig:UC_OPF}, exemplifies the optimal unit commitment plan for operating a notional power grid for the day-ahead, based on the available forecast of the demand and renewable generation. A cheap/slow unit provides the base generation for the entire planning horizon, while an expensive/fast generating unit is committed during the peak hours to avoid the violation of transmission limits.

The contributions of this paper are twofold. First the unit commitment problem is convexified via a family of linear and third-order semidefinite programming (TSDP) constraints. This convex relaxation achieves near-globally optimal solutions for UC problems with nearly $100,000$ binary variables. Second, a family of TSDP constraints are introduced to relax the power flow equations. The combination results in a tractable method for solving coupled UC-OPF optimization problems. The proposed method offers unprecedented scalability and improves upon the existing literature, in terms of the practical feasibility and efficiency of solutions, by allowing the joint optimization of commitment and power flow decisions based on an accurate nonlinear model for power grids. 



\subsection{Semidefinite Programming Relaxation}
Since a wide range of physical phenomena and dynamical systems can be modeled by polynomial functions, polynomial optimization have received significant research interest.
Our methodology is aligned with a popular framework for the study of polynomial optimization, that involves convex hull characterization of algebraic varieties through hierarchies of semidefinite programs (SDPs) \cite{lasserre2001global,lasserre2006convergent}.
Performance guarantees and extensions of SDP hierarchies have since been investigated by several papers \cite{laurent2009sums,gouveia2010theta,sojoudi2014exactness,belotti2015conic}, 
as well as their applications in various areas such as, quantum information theory \cite{tomamichel2016quantum,nagy2016epr}, compressed sensing \cite{javanmard2016phase,candes2015phase}, graph theory \cite{bandeira2016inference,aflalo2015convex}, statistics \cite{chandrasekaran2013computational}, operation of infrastructure networks \cite{lavaei2012zero,coogan2017offset}, and other branches of optimization theory \cite{park2015semidefinite}. The primary challenge for the application of SDP hierarchies beyond small-scale instances is the rapid growth of dimension, which necessitates a detailed study of sparsity and structural patterns to boost the efficiency 
\cite{muramatsu2003new,kim2003second,kim2003exact,bao2011semidefinite,natarajan2013penalized}.

\begin{table*}[!hp]
	\setlength{\tabcolsep}{2pt}
	\centering \small
		\caption{The performance of the TSDP relaxation algorithm for 24-hour horizon scheduling of benchmark systems with hourly epochs using the linear DC and nonlinear AC models.}\label{table_res_AC}\vspace{-2.0mm}
		\scalebox{0.8}{
	\begin{tabular}{lr||rrr|rr||rrr}
		\hline
		& &  \multicolumn{5}{c||}{Linear DC Model}  &  \multicolumn{3}{c}{Nonlinear AC Model}  \\
		\hline
		Test Case\!\!  & Number & Ratio of TSDP & TSDP & TSDP      & CPLEX & CPLEX  & Ratio of TSDP & TSDP & TSDP     \\
		& of Units & Inexact Binaries    & Gap	 & Time  & Gap & Time	  & Inexact Binaries    & Gap  & Time \\
		\hline
		IEEE 118\!\! 				& $54$ &  $0~/~1,296$  & $0\%$ & $3\mathrm{s}$	&	$0\%$	& $28\mathrm{s}$	& $0~/~1,296$ 		& $0.01\%$ & $11\mathrm{s}$\\
		IEEE 300\!\! 				& $69$ &  $0~/~1,656$  & $0\%$ & $3\mathrm{s}$	&	$0\%$	& $67\mathrm{s}$ & $0~/~1,656$ 		& $0.34\%$ & $41\mathrm{s}$\\
		\hline
		PEGASE 1354\!\! 			& $193$ & $28.4~/~4,632$  & $0.09\%$ & $18\mathrm{s}$	 & $8.57\%$ & $10,800\mathrm{s}^\dagger$ & $26.0~/~4,632$ 	& $1.24\%$ & $486\mathrm{s}$ \\
		PEGASE 2869\!\! 			& $392$ & $24.5~/~9,408$  & $0.01\%$ & $35\mathrm{s}$	 & $-$ & $10,800\mathrm{s}^\dagger$ & $33.8~/~9,408$ 	& $0.42\%$ & $2,175\mathrm{s}$ \\
		PEGASE 9241\!\! 			& $1,153$ & $75.4~/~27,672$ & $0.05\%$ & $137\mathrm{s}$	& $-$ & $10,800\mathrm{s}^\dagger$ & $226.5~/~27,672$ 	& $2.73\%$ & $56,351\mathrm{s}$ \\
		PEGASE 13659\!\!  			& $4,077$ & $29.5~/~97,848$ & $0.22\%$ & $266\mathrm{s}$ & $-$ & $10,800\mathrm{s}^\dagger$	& $995.3~/~97,848$ & $1.21\%$ & $116,064\mathrm{s}$ \\
		\hline
		\multicolumn{9}{l}{}\vspace{-6.5mm}\\ \\
		\multicolumn{9}{l}{\scriptsize $~^{\dagger}$ Solver is terminated after $3$ hours.}
	\end{tabular}
}
\end{table*}

\subsection{Review of Unit Commitment}
Economic scheduling of power generation units has been extensively investigated since the early 1960s, to handle predictable demand variations throughout the day-ahead. Extensions of the problem have later been studied to capture practical limits of network and security requirements, among other considerations. The reader is referred to \cite{allen2012price} for a detailed survey of the conventional formulations and computational methods for unit commitment. 

Recent policy and modeling proposals for electricity market operation and unit commitment include stochastic and robust optimization frameworks, under load and renewable generation uncertainty \cite{bitar2012bringing,bertsimas2013adaptive,phan2013optimization,yu2015impacts,lorca2016multistage,zhao2016unit,sundar2017modified}.
Additionally, incorporating other operational decisions into a comprehensive UC problem has been envisioned with a goal of co-optimizing multiple aspects of day-ahead planning, such as the optimal power flow \cite{bai2009semi,castillo2016unit,lipka2017running}, network topology control \cite{hedman2010co}, demand response \cite{wu2013hourly}, air quality control \cite{kerl2015new}, as well as scheduling of deferrable loads \cite{subramanian2013real}.

From a computational perspective, unit commitment algorithms rely on bounds from polynomial-time solvable relaxations for pruning search trees and certifying closeness to global optimality. Such relaxations can be generated through partial characterization of the convex hull of the feasible solutions \cite{ostrowski2012tight,lee2004min,damci2016polyhedral,geng2017alternative}. Additionally, in the presence of nonlinear price functions, conic inequalities can be adopted to strengthen the convex relaxations \cite{akturk2009strong,frangioni2009computational,bai2009semi,jabr2013rank}. Recently, a strong convex relaxation is proposed in \cite{fattahi2016conic} through a combination of reformulation-linearization and semidefinite programming techniques, which works very well on small instances of the unit commitment problem. Distributed methods are investigated in \cite{kargarian2015distributed} and \cite{papavasiliou2015applying} with the aim of leveraging high-performance computing platforms for solving large-scale unit commitment problems. Nevertheless, the improvements in run-time are reported to diminish with more than $15$ parallel workers \cite{papavasiliou2013comparative}. In terms of scalability, the proposed approach here significantly improves upon the above-referenced computational methods in the number of generating units as well as network size; notwithstanding, that our numerical experiments are conducted on a workstation with a single CPU.

\subsection{Review of Optimal Power Flow}

The optimal power flow problem is concerned with the determination of power flows and injections across the grid, for the optimal transmission and distribution of electricity. An accurate formulation of power flow in a transmission line includes nonconvex nonlinear equations, that substantively increase the computational complexity of the optimization problem. Consequently, the development of a framework for the joint optimization of UC and OPF has remained an open problem with significant economic impact as highlighted in \cite{baldick2005design}.

To this end, one of the most promising directions is based on the semidefinite programming relaxation of the power flow equations \cite{lavaei2012zero}. This approach to OPF has since been widely investigated and improved upon, through geometric analysis of feasible regions \cite{zhang2013geometry,kocuk2016inexactness,coffrin2016strengthening,chen2016bound,chen2017spatial}, and under certain graph-theoretic assumptions \cite{madani2015convex,gan2015exact}. Various studies have leveraged the sparsity of power networks for reducing the computational burden of solving semidefinite relaxations and developing distributed frameworks
\cite{molzahn2013implementation,andersen2014reduced,bose2015equivalent,madani2016promises,guo2016case,zhang2017distributed}.
More recently, other approaches such as Homotopy continuation \cite{mehta2016numerical}, for finding all solutions to power flow equations, and dynamic programming, in the presence of discrete variables \cite{dvijotham2017graphical} have been studied. Additionally, several extensions of OPF have been recently studied under more general settings, to address considerations such as the security of operation \cite{madani2016promises}, robustness \cite{dorfler2016breaking}, energy storage \cite{marley2016solving}, uncertainty of generation \cite{anese2017chance}. The reader is refereed to \cite{capitanescu2016critical} for a detailed survey on OPF.

\section*{Experimental Results}

This section gives a brief summary of the experiments with 
the proposed third-order semidefinite programming (TSDP) approach on large-scale instances of day-ahead scheduling. 
The goal is to determine the least-cost dispatch, that is, the on/off status and the amount of power produced by the generating units throughout the day ahead for meeting the load (demand) subject to the network transmission and technological constraints.
We consider 
real-world benchmark grids based on IEEE and European data with up to $13,659$ buses (vertices) and $4,077$ generating units. The 
planning horizon is divided into $24$ hourly intervals. For the largest benchmark, the model includes almost 100,000 binary decision variables.
Table \ref{table_res_AC} presents the average results for ten Monte Carlo demand simulations for each benchmark network.
The computations are performed on a workstation with a single CPU.
The details of data generation and experiments are discussed in the Methods and Materials section.

\subsection{Linear DC Model}
We first consider the approximate linear DC model, which is typically used by the electric power industry to formulate transmission of power in day-ahead scheduling problems.  
For all experiments, the proposed TSDP relaxation yields integer values for more than $99.5\%$ of binary variables. Moreover, the objective values of the recovered (feasible) scheduling decisions are provably within $0.22\%$ of global optimality for all benchmarks. The average performance of the TSDP relaxation, based on the DC model, is reported in columns three, four and five of Table \ref{table_res_AC}. Even for the largest benchmark, near-optimal solutions are obtained within a few minutes.

For comparison, the results with the commercial mixed-integer solver CPLEX, which is widely used by the system operators, are provided in columns six and seven. Although small-scale problems, based on IEEE data, are solved fast by CPLEX, no feasible solution is found after three hours of computation for the largest three benchmarks. 

\subsection{Nonlinear AC Model}
If an accurate nonlinear AC model for the flow of electricity is adopted, CPLEX, Gurobi and other commonly used off-the-shelf optimizers cannot be employed due to the presence of non-convex power constraints. For the largest benchmark system in Table \ref{table_res_AC}, the aforementioned nonlinear model results in a mixed-integer nonlinear optimization problem with $97,848$ binary variables, as well as $983,448$ non-convex quadratic constraints. For all experiments based on this network, our algorithm has been able to find solutions (with maximum power mismatch within $10^{-5}$  per-unit) that are on average within $2.73\%$ from global optimality. Moreover, for small- to medium-sized cases, all solutions are obtained in less than $40$ minutes and within $1.24\%$ gap from global optimality.

\section*{Notation}
The following notation is used in this paper: Bold letters are used for vectors and matrices, while italic letters with subscript indices refer to the entries of a vector or matrix. $\mathbb{R}$, $\mathbb{C}$, and $\mathbb{H}_n$ denote the sets of real numbers, complex numbers, and $n\times n$ Hermitian matrices, respectively. The letter ``$i$'' is reserved for the imaginary unit. The superscripts $(\cdot)^\ast$ and $(\cdot)^\top$ represent the conjugate transpose and transpose operators, respectively. The notations $\mathrm{real}\{\cdot\}$, $\mathrm{imag}\{\cdot\}$, and $|\cdot|$ represent the real part, imaginary part, and element-wise absolute value of a scalar or matrix, respectively. The notations $\mathbf{X}_{\bullet,k}$ and $\mathbf{X}_{l,\bullet}$ refer to the $k$-th column and the $l$-th row of matrix $\mathbf{X}$, respectively. Additionally, $\mathrm{diag}\{\mathbf{X}\}$ denotes the vertical vector whose entries are given by the diagonal elements of $\mathbf{X}$. The notation $\mathbf{X}\succeq0$ means that $\mathbf{X}$ is Hermitian and positive semidefinite. The notation $\mathbf{X}\otimes\mathbf{Y}$ refers to the Kronecker product of the matrices $\mathbf{X}$ and $\mathbf{Y}$. 
Given an $n\times n$ matrix $\mathbf{X}$ and $\mathcal{S}_1,\mathcal{S}_2\in\{1,\ldots,n\}$, define $\mathbf{X}[\mathcal{S}_1,\mathcal{S}_2]$ to be the $|\mathcal{S}_1|\times|\mathcal{S}_2|$ submatrix of $\mathbf{X}$ with row and columns from $\mathcal{S}_1$ and $\mathcal{S}_2$, respectively. Throughout the paper, non-positive subscript indices refer to known initial values.

\section*{Power System Scheduling}
The power system scheduling problem seeks to find the most economic operation plan for a set of generating units throughout a time horizon to meet the demand for electricity, subject to technological constraints. Let $\mathcal{G}=\{1,2,\ldots,G\}$ denote the set of generating units, whose schedule and the amount of contribution to the grid are to be determined. In order to formulate the problem as a static optimization, it is common practice to divide the planning horizon into a set of discrete time intervals $\mathcal{T}=\{1,2,\ldots,T\}$, e.g., hourly time slots for day-ahead scheduling.

Let $x_{g,t}\in\{0,1\}$ be a binary variable indicating the decision of whether or not the generating unit $g\in\mathcal{G}$ is committed for production in the time slot $t\in\mathcal{T}$. If $x_{g,t}=1$, the unit is active and expected to produce power within its capacity limitations, otherwise, no power is produced by $g$ in this time slot. Additionally, let $c_{g,t}$ be the production cost of unit $g$, during the interval $t$. There are two types of power exchanges between generating units and loads in a power system: i) active power, and ii) reactive power. Active power is the actual product that is traded to meet the demand, whereas the reactive power is a technical term, which represents the oscillation exchanges between generators and loads that help maintaining voltages. Let $p_{g,t}$ and $q_{g,t}$, respectively, to be the amount of active power and reactive powers produced by unit $g\in\mathcal{G}$, in time interval $t\in\mathcal{T}$. The overall power injection by generating unit $g\in\mathcal{G}$ can be expressed as the complex number $p_{g,t} + i q_{g,t}$, which is referred to as {\it complex power}, where ``$i$'' accounts for the imaginary unit.

A distinctive feature of our approach is the ability to jointly optimize unit commitment and power flow for an accurate model of the grid. In this paper, the constraints of this large-scale optimization problem are divided into two classes: 
\begin{enumerate}
	\item {\bf Unit constraints} that model the capacity and technological limitations of generating units, and 
	\item {\bf Network constraints} that model laws of physics governing the flow of power electricity across the grid, such as conservation of energy, as well as the transmission capacity limitations and demand requirements, throughout the planning horizon.
\end{enumerate}

Using the notation introduced above, we formulate power system scheduling as the following optimization problem: 
\begin{subequations}\label{PSS}\begin{align}
		& \underset{
			\begin{subarray}{c}
				\\
				\mathbf{x},\mathbf{p},\mathbf{q},\mathbf{c}\in\mathbb{R}^{G\! \times\! T}\!\!\!\!\!\!\!\!\!\! 
			\end{subarray}
		}{\text{minimize}}
		& & \sum_{g=1}^{G}{\sum_{t=1}^{T}{c_{g,t}}} &&&&  \label{PSS_obj}\\
		& \text{subject to}
		& & (\mathbf{x}^{\top}_{g,\bullet},\mathbf{p}^{\top}_{g,\bullet},\mathbf{q}^{\top}_{g,\bullet},\mathbf{c}^{\top}_{g,\bullet})\in\mathcal{U}_g &&&& \hspace{-0.4cm}\forall g\in\mathcal{G}, \label{cons_unit}\\
		& & & (\mathbf{p}_{\bullet, t},\mathbf{q}_{\bullet, t})\in\mathcal{N}_t &&&& \hspace{-0.36cm}\forall t\in\mathcal{T},\label{cons_network}
\end{align}\end{subequations}
with respect to decision variables
$\mathbf{x}\equiv [x_{g,t}]$,
$\mathbf{c}\equiv [c_{g,t}]$,
$\mathbf{p}\equiv[p_{g,t}]$ and 
$\mathbf{q}\equiv[q_{g,t}]$. Optimization problem [\ref{PSS_obj}]--[\ref{cons_network}] minimizes the overall cost of producing power subject to unit and network constraints [\ref{cons_unit}] and [\ref{cons_network}], respectively. For every generating unit $g\in\mathcal{G}$, the quadruplet $(\mathbf{x}^{\top}_{g,\bullet},\mathbf{p}^{\top}_{g,\bullet},\mathbf{q}^{\top}_{g,\bullet},\mathbf{c}^{\top}_{g,\bullet})\in\mathbb{R}^{T\times 4}$ characterizes the scheduling decision, throughout the planning horizon, while for every time slot $t\in\mathcal{T}$, the pair $(\mathbf{p}_{\bullet, t},\mathbf{q}_{\bullet, t})\in\mathbb{R}^{G\times 2}$ accounts for the generation profile. The price functions and technological limitations of generating units are described by the sets $\mathcal{U}_1,\mathcal{U}_2,\ldots,\mathcal{U}_G\subset\mathbb{R}^{T\times 4}$, while the demand information and network data across the time slots are given by $\mathcal{N}_1,\mathcal{N}_2,\ldots,\mathcal{N}_T\subset\mathbb{R}^{G\times 2}$.





The binary unit commitment decisions and nonlinearity of network equations are the primary sources of computational complexity for solving the problem [\ref{PSS_obj}]--[\ref{cons_network}]. 
As a result, there has been a huge body of research devoted to finding convex relaxations for power system scheduling and its related problems, by means of tools and techniques from the area of mathematical programming. In the following, we first describe the families of sets $\{\mathcal{U}_g\}_{g\in\mathcal{G}}$ and $\{\mathcal{N}_t\}_{t\in\mathcal{T}}$, given by the unit commitment and network constraints of power system scheduling. We then introduce convex surrogates for them which lead to a class of computationally tractable and, yet, accurate relaxations of problem [\ref{PSS_obj}]--[\ref{cons_network}]. 

\subsection{Unit Constraints}
Following is a definition for the family $\{\mathcal{U}_g\}_{g\in\mathcal{G}}$, which is based on a number of practical limitation for the operation of generating units.
\begin{definition}
For every generating unit $g\in\mathcal{G}$, define $\mathcal{U}_g$ to be the set of all quadruplets
$(\mathbf{x}^{\top}_{g,\bullet},\mathbf{p}^{\top}_{g,\bullet},\mathbf{q}^{\top}_{g,\bullet},\mathbf{c}^{\top}_{g,\bullet})\in\mathbb{R}^{T\times 4}$ that satisfies constraints {\rm[\ref{cons_bin}]}, {\rm[\ref{cons_cost}]}, {\rm[\ref{cons_cap}]}, {\rm[\ref{cons_min}]}, and {\rm[\ref{cons_ramp}]}, for all $t\in\mathcal{T}$.
\end{definition}
Note that, non-positive indices refer to given initial values.
In the reminder of this section, we detail each of the above-mentioned constraints.


\paragraph{Production Costs:}
The cost of operating a unit $g\in\mathcal{G}$ within different time intervals is a quadratic function of the active power produced by the unit. In addition, there is a fixed cost $\gamma_g$ associated with every interval during which the generator is committed (i.e., $x_{g,t}=1$), as well as a startup cost $\gamma_g^{\uparrow}$ and a shutdown cost $\gamma_g^{\downarrow}$ that are enforced on time slots at which the unit $g$ changes status. Therefore, the price of operating unit $g$ at time $t$ can be described through the nonlinear equation [\ref{cons_cost}],
where $\alpha_g$ and $\beta_g$ are nonnegative coefficients. 

\paragraph{Generation Capacity:}
If a generating unit $g\in\mathcal{G}$ is committed at time $t\in\mathcal{T}$, the amount of active power $p_{g,t}$ and reactive power $q_{g,t}$ produced in that time slot must lie within capacity limitations of the unit. In other words, if $x_{g,t}=1$, then we have $p_{g,t}\in[\ubar{p}_g,\bar{p}_g]$ and $q_{g,t}\in[\ubar{q}_g,\bar{q}_g]$, where $\ubar{p}_g$, $\bar{p}_g$, $\ubar{q}_g$ and $\bar{q}_g$ are the given lower and upper bounds for unit $g$. 
Constraints [\ref{cons_cap_p}]--[\ref{cons_cap_q}], ensure that the amount of power produced by unit $g$ is zero if $x_{g,t}=0$, and within capacity limits, if $x_{g,t}=1$.

\paragraph{Minimum Up \& Down Time Limits:}
Technical considerations often prohibit frequent changes in the status of generating units. Once a unit starts producing power, there is a minimum time before it can be turned off, and once the unit is turned off, it cannot be immediately activated, again. Denote by $m^{\uparrow}_g$ and $m^{\downarrow}_g$ the minimum time for which the generating unit $g\in\mathcal{G}$ is required to remain active and deactivate, respectively. The minimum up and down limits for unit $g$ are enforced through constraints [\ref{cons_min_u}]--[\ref{cons_min_l}].

\paragraph{Ramp Rate Limits:}
The rate of change in the amount of power produced by a generating unit is often constrained, depending on the type of the generator. Denote by $r_g$ the maximum variation of active power generation, that is allowed by unit $g\in\mathcal{G}$, between two adjacent time intervals in which the unit is committed. Similarly, define $s_g$ as the maximum amount of active power that can be generated by unit $g$ immediately after startup or prior to shutdown. Ramp rate limits of unit $g\in\mathcal{G}$ are expressed through constraints [\ref{cons_ramp_pu}] and [\ref{cons_ramp_pl}].
Observe that if either $x_{g,t-1} = 0$ or $x_{g,t} = 0$, the constraints in [\ref{cons_ramp_pu}] and [\ref{cons_ramp_pl}] reduce to
$|p_{g,t}|\leq s_g$. Alternatively, if $x_{g,t-1} = x_{g,t} = 1$, the above constraints imply that $|p_{g,t}-p_{g,t-1}|\leq r_g$.


\begin{table}[t]
	\vspace{1.3mm}
	\line(1,0){360}\\
	\begin{flushleft}
	{\bf Unit Constraints:}
		\end{flushleft}
	\begin{align}
	x_{g,t} \in \{0,1\}\label{cons_bin}
	\end{align}\vspace{-0.5cm}
	\begin{align}	
	&c_{g,t} = \alpha_g p^2_{g,t} + \beta_g p_{g,t} +\nonumber\\ 
	&\quad\quad\;\;\,\gamma_g x_{g,t} + \gamma_g^{\uparrow} (1-x_{g,{t-1}})x_{g,t} + \gamma_g^{\downarrow}	 x_{g,t-1}(1-x_{g,t}),\label{cons_cost}
	\end{align}\vspace{-0.5cm}
	\begin{subequations}\label{cons_cap}\begin{align}	
		& \ubar{p}_g x_{g,t} \leq  p_{g,t} \leq \bar{p}_g x_{g,t}\label{cons_cap_p}\\ 
		& \ubar{q}_g x_{g,t} \leq  q_{g,t} \leq \bar{q}_g x_{g,t}.\label{cons_cap_q}
		\end{align}\end{subequations}\vspace{-0.5cm}
	\begin{subequations}\label{cons_min}\begin{align}	
		&\quad\;\;\, x_{g,t}   \geq x_{g,\tau}-x_{g,\tau-1},  \quad \forall\tau\in\{t-m^{\uparrow}_g+1,\ldots,t\},\label{cons_min_u}\\
		&1-x_{g,t} \geq x_{g,\tau-1}-x_{g,\tau},  \quad \forall\tau\in\{t-m^{\downarrow}_g+1,\ldots,t\}.\label{cons_min_l}
		\end{align}\end{subequations}\vspace{-0.5cm}
	\begin{subequations}\label{cons_ramp}\begin{align}		
		&p_{g,t} - p_{g,t-1} \leq  r_g x_{g,t-1} + s_g (1-x_{g,t-1}),\label{cons_ramp_pu} \\ 
		&p_{g,t-1} - p_{g,t} \leq r_g x_{g,t} + s_g (1-x_{g,t}), \label{cons_ramp_pl}
		\end{align}\end{subequations}

	\vspace{-3mm}	
	
	\line(1,0){360}\\
\begin{flushleft}
	{\bf AC Network Constraints:}
\end{flushleft}
	\begin{subequations}\label{cons_AC}\begin{align}
		\mathbf{d}_{t}+\mathrm{diag}\{\mathbf{v}_{\bullet, t}\mathbf{v}^{\ast}_{\bullet, t}\mathbf{Y}^{\ast}_t\}&=\mathbf{C}^{\top}(\mathbf{p}_{\bullet, t}+i\mathbf{q}_{\bullet, t})\label{cons_AC1}\\
		\lvert\mathrm{diag}\{\vec{\mathbf{C}}_t\mathbf{v}_{\bullet, t}\mathbf{v}^{\ast}_{\bullet, t}\vec{\mathbf{Y}}_t^{\ast}\}\rvert&\leq \mathbf{f}_{\mathrm{max};t}\label{cons_AC2}\\
		\lvert\mathrm{diag}\{\cev{\mathbf{C}}_t\mathbf{v}_{\bullet, t}\mathbf{v}^{\ast}_{\bullet, t}\cev{\mathbf{Y}}_t^{\ast}\}\rvert&\leq \mathbf{f}_{\mathrm{max};t}\label{cons_AC3}\\
		\mathbf{v}_{\mathrm{min}}\leq\lvert\mathbf{v}_{\bullet, t}\rvert&\leq \mathbf{v}_{\mathrm{max}}\label{cons_AC4}
		\end{align}\end{subequations}
	
	\vspace{-3mm}
	
	\line(1,0){360}\\
\begin{flushleft}
	{\bf DC Network Constraints:}
\end{flushleft}
	\begin{subequations}\label{cons_DC}\begin{align}
		\mathbf{q}_{\bullet, t}&=0\label{cons_DC1}\\
		\mathrm{real}\{\mathbf{d}_{t}\}+\mathbf{B}_t \boldsymbol{\theta}_{\bullet, t} &=\mathbf{C}^{\top} \mathbf{p}_{\bullet, t}\label{cons_DC2}\\
		\lvert \vec{\mathbf{B}}_t \boldsymbol{\theta}_{\bullet, t} \rvert&\leq \mathbf{f}_{\mathrm{max};t}\label{cons_DC3}
		\end{align}\end{subequations}
	
	\vspace{-3mm}
	
	\line(1,0){360}
	\caption{Unit and network constraints in power system scheduling.}
	\label{fig:table_cons}
	
	\vspace{-5mm}
\end{table}

\subsection{Network Constraints}
In this part, we focus on network considerations in power system scheduling. The transmission of electricity from suppliers to consumers is carried out through an interconnected network whose topology throughout each time interval $t\in\mathcal{T}$ can be modeled as a directed graph $\mathcal{H}_t=(\mathcal{V},\mathcal{E}_t)$, with $\mathcal{V}$ and $\mathcal{E}_t$ as the set of vertices and edges, respectively. In power system terminology, vertices are referred to as ``buses'', and edges are called ``lines'' or ``branches'' of the network. Each generating unit is associated with (located at) one of the buses. Define the unit incidence matrix $\mathbf{C}\in\{0,1\}^{G\times\mathcal{V}}$ to be a binary matrix whose entry $(g,k)$ is equal to one, if and only if the generating unit $g$ belongs to bus $k$. Additionally, define the pair of matrices $\vec{\mathbf{C}}_t,\cev{\mathbf{C}}_t\in\{0,1\}^{\mathcal{E}_t\times\mathcal{V}}$ as the initial and final incidence matrices, respectively. The entry $(k,l)$ of $\vec{\mathbf{C}}_t$ is equal to one, if and only if line $l$ starts at bus $k$, while the entry $(k,l)$ of $\cev{\mathbf{C}}_t$ equals one, if and only if line $l$ ends at bus $k$. 

The steady state voltages across the network are sinusoidal functions with a global frequency. As a result, the voltage function at each bus can be characterized by its amplitude and phase difference from a reference bus. Therefore, for each $k\in\mathcal{V}$ and $t\in\mathcal{T}$, a complex number $v_{k;t}$ is defined, whose magnitude $|v_{k;t}|$ and angle $\angle v_{k;t}$, respectively, account for the amplitude and phase of the voltage at the bus $k$, in time interval $t$. 
Define $\mathbf{v}\triangleq[v_{k;t}]\in\mathbb{C}^{\mathcal{V}\times \mathcal{T}}$ and $\boldsymbol{\theta}\triangleq[\angle v_{k;t}]\in\mathbb{R}^{\mathcal{V}\times \mathcal{T}}$, to be the matrices encapsulating complex voltage and phase angle values, respectively.

\begin{figure}
	\centering
	\includegraphics[width=0.75\linewidth]{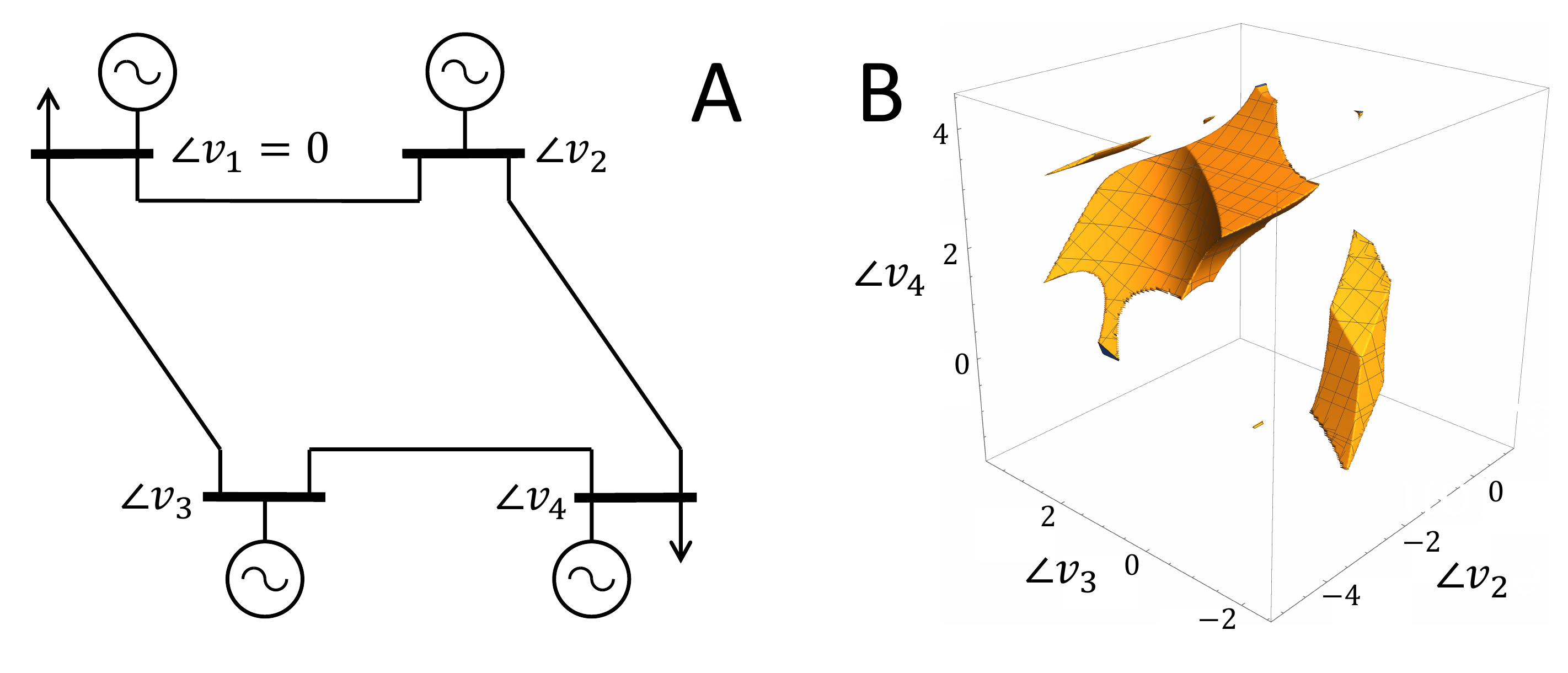}
	\caption{ (A) A four bus power system from \cite{zimmerman2011matpower} with two loads and four generators. (B) Feasible region of voltage angles, in which demand and technological constraints of AC formulation are satisfied.}
	\label{fig:FigureAC}
\end{figure}

The two widely used models for power networks are discussed next. The first one is the accurate Alternating Current (AC) model, which incorporates the nonlinear power flow equations. The next one is the Direct Current (DC) model, which is a simplified version of the AC model and can be described by linear equalities. 
The use of nonlinear AC power flow equations introduces substantial complexity into power system optimization problems. However, various physical phenomena, such as network losses and reactive power flows are captured by the AC model, while ignored by the DC model. As a result, it is desirable to adopt the AC model, in order to determine better operation strategies. 
Figure \ref{fig:FigureAC} illustrates a highly non-convex feasible region of voltage angles, enforced by the demand and technological constraints, in a simple four bus network that is described exactly by the AC model.
One of the primary benefits of the proposed method in this paper, is the possibility of adopting the AC model in large-scale power system scheduling problems.

\subsubsection{Alternating Current Power Flow Model}

In the AC model, characteristics of the network in a time interval $t\in\mathcal{T}$, can be described by a triplet of admittance matrices $\vec{\mathbf{Y}}_t,\cev{\mathbf{Y}}_t\in\mathbb{C}^{\mathcal{E}_t\times \mathcal{V}}$ and $\mathbf{Y}_t\in\mathbb{C}^{\mathcal{V}\times \mathcal{V}}$, that govern the flow of power throughout the network. 
Next, we define the family $\{\mathcal{N}^{\mathrm{AC}}_t\}_{t\in\mathcal{T}}$
and give a brief description for each constraint.
\begin{definition}
	For every time interval $t\in\mathcal{T}$, let $\mathcal{N}^{\mathrm{AC}}_t$ be the set of pairs
	$(\mathbf{p}_{\bullet,t},\mathbf{q}_{\bullet,t})\in\mathbb{R}^{T\times 2}$, for which there exists a vector of complex voltages $\mathbf{v}_{\bullet,t}\in\mathbb{C}^{\mathcal{V}}$ satisfying the constraints {\rm[\ref{cons_AC1}]}--{\rm[\ref{cons_AC4}]}.
\end{definition}
\paragraph{AC Power Balance Equation:}
Constraint [\ref{cons_AC1}] is referred to as the {\it power balance equation} which accounts for the conservation of energy at all buses of the network. The vector $\mathbf{d}_t\in\mathbb{C}^{\mathcal{N}}$ denotes the demand forecast at each bus, in interval $t$, whose real and imaginary parts account for active and reactive power demands, respectively. Observe that the overall complex power produced by generating units located at each bus $k\in\mathcal{V}$ is given by the $k$-th entry of $\mathbf{C}^{\top}(\mathbf{p}_{\bullet t}+i\mathbf{q}_{\bullet t})$. Finally, the $k$-th entry of the vector $\mathrm{diag}\{\mathbf{v}_{\bullet, t}\mathbf{v}^{\ast}_{\bullet, t}\mathbf{Y}^{\ast}_t\}$ is equal to the amount of complex power exchange between bus $k$ and the rest of the network. The voltages across the network are adjusted in such a way that the overall complex power produced at each bus equals the sum of power consumptions and power exchanges of that bus, at all times. This requirement is enforces by constraint [\ref{cons_AC1}].

\paragraph{AC Thermal Limits:}
Due to thermal losses, the flow entering a line may differ from the flow leaving the line at the other end. For each time interval $t\in\mathcal{T}$, complex power flows entering the lines of the network through their starting and ending buses are given by vectors $\mathrm{diag}\{\vec{\mathbf{C}}_t\mathbf{v}_{\bullet, t}\mathbf{v}^{\ast}_{\bullet, t}\vec{\mathbf{Y}}_t^{\ast}\}$ and $\mathrm{diag}\{\cev{\mathbf{C}}_t\mathbf{v}_{\bullet, t}\mathbf{v}^{\ast}_{\bullet, t}\cev{\mathbf{Y}}_t^{\ast}\}$, respectively. Constraints [\ref{cons_AC2}] and [\ref{cons_AC3}] restrict the flow of power, within the thermal limit of the lines $\mathbf{f}_{\mathrm{max};t}\in\mathbb{R}^{\mathcal{L}_t}$, for each $t\in\mathcal{T}$.

\paragraph{Voltage Magnitude Limits:}
In order for power system components to operate properly, the voltage magnitude at each bus needs to remain within a prespecified range, given by vectors $\mathbf{v}_{\mathrm{min}},\mathbf{v}_{\mathrm{max}}\in\mathbb{R}^{\mathcal{V}}$. Voltage magnitude limits are enforced through the constraint [\ref{cons_AC4}].

The nonlinear AC network constraints [\ref{cons_AC1}]--[\ref{cons_AC4}] pose a significant challenge for solving power system optimization problems based on a full model. As a result, typically a simplified version of the AC model is considered in practice which is explained next.

\subsubsection{Direct Current Power Flow Model}
The DC model can be formulated by ignoring the reactive powers, voltage magnitude deviations from their nominal values, and network losses. Under this model, the network is described by means of sustenance matrices $\vec{\mathbf{B}}_t\in\mathbb{R}^{\mathcal{E}_t\times \mathcal{V}}$ and $\mathbf{B}_t\in\mathbb{R}^{\mathcal{V}\times \mathcal{V}}$. Moreover, the flow of active power across the network in each time interval $t\in\mathcal{T}$ is expressed with respect to the vector of voltage angles $\boldsymbol{\theta}_{\bullet,t}\in\mathbb{R}^{\mathcal{V}}$.
\begin{definition}
	For every time interval $t\in\mathcal{T}$, let $\mathcal{N}^{\mathrm{DC}}_t$ be the set of pairs
	$(\mathbf{p}_{\bullet,t},\mathbf{q}_{\bullet,t})\in\mathbb{R}^{T\times 2}$, for which there exists a vector of voltage phase values $\boldsymbol{\theta}_{\bullet,t}\in\mathbb{R}^{\mathcal{V}}$ satisfying constraints {\rm[\ref{cons_DC1}]}--{\rm[\ref{cons_DC3}]}.
\end{definition}
Constraint [\ref{cons_DC2}] is a simplified alternative for power balance equation [\ref{cons_AC1}], in which the vector $\mathbf{B}_t \boldsymbol{\theta}_{\bullet, t}\in\mathbb{R}^{\mathcal{V}}$ contains approximate values for active power exchanges between each vertex and the rest of the network. Additionally, thermal limits are enforced through the constraint [\ref{cons_DC3}], in which $\vec{\mathbf{B}}_t \boldsymbol{\theta}_{\bullet, t}\in\mathbb{R}^{\mathcal{L}_t}$ is the vector of approximate values for active power flow of lines. Notice that, since the network losses are ignored under the DC model, power flows entering both directions are considered equal and it suffices to enforce one inequality for each line.



\section*{Methodology}
In order to tackle a general power system scheduling problem of the form [\ref{PSS_obj}]--[\ref{cons_network}], we develop third-order semidefinite programming (TSDP) relaxations for the families of sets $\{\mathcal{U}_g\}_{g\in\mathcal{G}}$ and $\{\mathcal{N}_t\}_{t\in\mathcal{T}}$, which lead to a computationally-tractable algorithm. The proposed approach involves introducing additional variables, each as a proxy for a quadratic monomial. We design a class of inequalities, to strengthening the relation between each proxy variable and the monomial it represents.

\begin{table*}[t]
			\vspace{1mm}\line(1,0){360}
	\begin{subequations}\label{cons_KK}
		\begin{align}\label{cons_K}
		\!\!\!\!\!\!\!
		\left(\!\!
		\begin{bmatrix}
		& \hspace{-1.7mm} +1 &\hspace{-1.7mm}  \\
		1 &\hspace{-1.7mm} -1 &\hspace{-1.7mm} \\	
		&\hspace{-1.7mm}  &\hspace{-1.7mm} +1  \\
		1 &\hspace{-1.7mm}  &\hspace{-1.7mm} -1							
		\end{bmatrix}\!\!\otimes\!\!
		\begin{bmatrix}
		&\hspace{-1mm} -\ubar{p}_g 			  	&\hspace{-1.2mm}  				&\hspace{-1.2mm} +1  &\hspace{-1mm}   \\
		&\hspace{-1mm} +\bar{p}_g 				  &\hspace{-1.2mm}  					&\hspace{-1.2mm} -1  &\hspace{-1mm}   \\
		&\hspace{-1mm}  				 				&\hspace{-1.2mm} -\ubar{p}_g &\hspace{-1.2mm}   &\hspace{-1mm} +1  \\
		&\hspace{-1mm}  				 				&\hspace{-1.2mm} +\bar{p}_g &\hspace{-1.2mm}   &\hspace{-1mm} -1  \\
		s_g&\hspace{-1.2mm} r_g\!-\!s_g &\hspace{-1mm}  			&\hspace{-1.2mm} +1  &\hspace{-1mm} -1  \\
		s_g&\hspace{-1.2mm}  								&\hspace{-1mm} r_g\!-\!s_g &\hspace{-1.2mm} -1  &\hspace{-1mm} +1 \\								
		\end{bmatrix}\right)
		\!\!\times\!\!
		\begin{bmatrix}
		\mathbf{e}^{\top}_1\\
		\mathbf{e}^{\top}_2\!\!+\!\mathbf{e}^{\top}_{6}\!\!+\!\mathbf{e}^{\top}_{7}\\
		\mathbf{e}^{\top}_3\!\!+\!\mathbf{e}^{\top}_{11}\!\!+\!\mathbf{e}^{\top}_{13}\\
		\mathbf{e}^{\top}_4\!+\!\mathbf{e}^{\top}_{9}\\
		\mathbf{e}^{\top}_5\!+\!\mathbf{e}^{\top}_{15}\\
		\mathbf{e}^{\top}_8\!+\!\mathbf{e}^{\top}_{12}\\
		\mathbf{e}^{\top}_{14}\\
		\mathbf{e}^{\top}_{10}
		\end{bmatrix}^{\!\!\top}\!\!\!\!\times\!\!
		\begin{bmatrix}
		1\\
		x_{g,t-1}\\
		x_{g,t}\\
		p_{g,t-1}\\
		p_{g,t}\\
		u_{g,t}\\
		y_{g,t}\\
		z_{g,t}
		\end{bmatrix}\!\!\geq\! 0,\!\!\!\!
		\end{align}
		\vspace{1mm}\line(1,0){360}
		\begin{align}\label{cons_Kd}
		\!\!\!\!\!\!\!
		\left(
		\begin{bmatrix}
		& \hspace{-1.5mm} +1 &\hspace{-1.5mm}  \\
		1 &\hspace{-1.5mm} -1 &\hspace{-1.5mm} \\	
		&\hspace{-1.5mm} +\bar{p}_g     &\hspace{-1.5mm} -1	\\
		&\hspace{-1.5mm} -\ubar{p}_g   &\hspace{-1.5mm} +1 						
		\end{bmatrix}\!\!\otimes\!\!
		\begin{bmatrix}
		+1\\
		-1\\
		+1\\
		-1
		\end{bmatrix}^{\!\!\top}
		\right)
		\!\!\times\!\!
		\begin{bmatrix}
		\dot{\mathbf{e}}^{\top}_1\\
		\dot{\mathbf{e}}^{\top}_2\\
		\dot{\mathbf{e}}^{\top}_3+\dot{\mathbf{e}}^{\top}_{5}+\dot{\mathbf{e}}^{\top}_{7}\!\!\\
		\dot{\mathbf{e}}^{\top}_4\\
		\dot{\mathbf{e}}^{\top}_6\\
		\dot{\mathbf{e}}^{\top}_8\\
		\dot{\mathbf{e}}^{\top}_{9}+\dot{\mathbf{e}}^{\top}_{11}\\
		\dot{\mathbf{e}}^{\top}_{10}\\
		\dot{\mathbf{e}}^{\top}_{12}
		\end{bmatrix}^{\top} \hspace{-2.0mm}\times
		\begin{bmatrix}
		1\\
		x_{g,t-2}\\
		x_{g,t-1}\\
		x_{g,t}\\
		u_{g,t-1}\\
		u_{g,t}\\
		p_{g,t-1}\\
		z_{g,t-1}\\
		y_{g,t}
		\end{bmatrix}
		\geq 0,
		\end{align}
		\vspace{1mm}\line(1,0){360}
		\begin{align}\label{cons_Ku}
		\!\!\!\!\!\!\!
		\left(
		\begin{bmatrix}
		& \hspace{-1.5mm} +1 &\hspace{-1.5mm}  \\
		1 &\hspace{-1.5mm} -1 &\hspace{-1.5mm} \\	
		&\hspace{-1.5mm} +\bar{p}_g     &\hspace{-1.5mm} -1	\\
		&\hspace{-1.5mm} -\ubar{p}_g   &\hspace{-1.5mm} +1 						
		\end{bmatrix}\!\!\otimes\!\!
		\begin{bmatrix}
		-1\\
		+1\\
		-1
		\end{bmatrix}^{\top}
		\right)
		\!\!\times\!\!
		\begin{bmatrix}
		\ddot{\mathbf{e}}^{\top}_1\\
		\ddot{\mathbf{e}}^{\top}_2+\ddot{\mathbf{e}}^{\top}_{5}\\
		\ddot{\mathbf{e}}^{\top}_3\\
		\ddot{\mathbf{e}}^{\top}_4\\
		\ddot{\mathbf{e}}^{\top}_6\\
		\ddot{\mathbf{e}}^{\top}_8\\
		\ddot{\mathbf{e}}^{\top}_7\\
		\ddot{\mathbf{e}}^{\top}_9
		\end{bmatrix}^{\top} \hspace{-2.0mm}\!\!\times\!
		\begin{bmatrix}
		x_{g,t-2}\\
		x_{g,t-1}\\
		x_{g,t}\\
		u_{g,t-1}\\
		u_{g,t}\\
		p_{g,t-1}\\
		z_{g,t-1}\\
		y_{g,t}
		\end{bmatrix}\geq 0.
		\end{align}
		\vspace{1mm}\line(1,0){360}\\
		\begin{flushleft}
		{\scriptsize $\otimes$ denotes the Kronecker product of two matrices. 
			$\left\{\mathbf{e}_1,\mathbf{e}_2,\ldots,\mathbf{e}_{15}\right\}$,
			$\left\{\dot{\mathbf{e}}_1,\dot{\mathbf{e}}_2,\ldots,\dot{\mathbf{e}}_{12}\right\}$ and
			$\left\{\ddot{\mathbf{e}}_1,\ddot{\mathbf{e}}_2,\ldots,\ddot{\mathbf{e}}_{9}\right\}$
			denote the standard basis vectors for 
			$\mathbb{R}^{15}$, $\mathbb{R}^{12}$, and $\mathbb{R}^{9}$,
			respectively. 
		}		\end{flushleft}
	\end{subequations}
\end{table*}

In this work, we propose a convex relaxation of the power system scheduling problem [\ref{PSS_obj}]--[\ref{cons_network}], which is built by substituting the unit and AC network feasible sets  with their convex surrogates  
$\{\mathcal{U}^{\mathrm{TSDP}}_g\}_{g\in\mathcal{G}}$ and 
$\{\mathcal{N}^{\mathrm{TSDP}}_t\}_{t\in\mathcal{T}}$, respectively:
\begin{subequations}\label{PSS_relax}\begin{align}
	& \underset{
		\begin{subarray}{c}
		\\
		\mathbf{x},\mathbf{p},\mathbf{q},\mathbf{c}\in\mathbb{R}^{G\! \times\! T}\!\!\!\!\!\!\!\!\!\! 
		\end{subarray}
	}{\text{minimize}}
	& & \sum_{g=1}^{G}{\sum_{t=1}^{T}{c_{g,t}}} &&&&  \label{PSS_relax_obj}\\
	& \text{subject to}
	& & (\mathbf{x}^{\top}_{g,\bullet},\mathbf{p}^{\top}_{g,\bullet},\mathbf{q}^{\top}_{g,\bullet},\mathbf{c}^{\top}_{g,\bullet})\in\mathcal{U}^{\mathrm{TSDP}}_g &&&& \hspace{-0.4cm}\forall g\in\mathcal{G}, \label{cons_unit_relax}\\
	& & & (\mathbf{p}_{\bullet, t},\mathbf{q}_{\bullet, t})\in\mathcal{N}^{\mathrm{TSDP}}_t &&&& \hspace{-0.36cm}\forall t\in\mathcal{T},\label{cons_network_relax}
\end{align}\end{subequations}
Due to convexity of the sets $\{\mathcal{U}^{\mathrm{TSDP}}_g\}_{g\in\mathcal{G}}$ and 
$\{\mathcal{N}^{\mathrm{TSDP}}_t\}_{t\in\mathcal{T}}$, the problem [\ref{PSS_relax_obj}]--[\ref{cons_network_relax}] can be solved in polynomial time. Moreover, since $\mathcal{U}_g\subseteq\mathcal{U}^{\mathrm{TSDP}}_g$ and $\mathcal{N}^{\mathrm{AC}}_t\subseteq\mathcal{N}^{\mathrm{TSDP}}_t$, for every $g\in\mathcal{G}$ and $t\in\mathcal{T}$, respectively, the optimal cost of problem [\ref{PSS_relax_obj}]--[\ref{cons_network_relax}] is a lower bound to the optimal cost of problem [\ref{PSS_obj}]--[\ref{cons_network}]. If an optimal solution to the problem [\ref{PSS_relax_obj}]--[\ref{cons_network_relax}] satisfies the original constraints [\ref{cons_unit}] and [\ref{cons_network}], then the relaxation is exact and a provably global optimal solution to problem [\ref{PSS_obj}]--[\ref{cons_network}] is obtained. Otherwise, a rounding procedure is adopted to transform the optimal solution of [\ref{PSS_relax_obj}]--[\ref{cons_network_relax}] to a feasible and near optimal solution of [\ref{PSS_obj}]--[\ref{cons_network}].

\subsection{Relaxation of Unit Constraints}
Each unit feasible set $\mathcal{U}_g$ is a semialgebraic set, with constraints [\ref{cons_bin}] and [\ref{cons_cost}] as the sources of nonconvexity. 
In this work, we create a family of convex surrogates $\{\mathcal{U}^{\mathrm{TSDP}}_g\}_{g\in\mathcal{G}}$, by enforcing a collection of linear and conic inequalities. 
To this end, define auxiliary variables
$\mathbf{u},
\mathbf{y},
\mathbf{z},
\mathbf{o}\in\mathbb{R}^{G\times T}$, whose components account for monomials $x_{g,t-1}x_{g,t}$,
$p_{g,t-1}x_{g,t}$,
$x_{g,t-1}p_{g,t}$ and
$p^2_{g,t}$, respectively. In other words, if the relaxation is exact, the equations
\begin{subequations}
\begin{align}
&u_{g,t}=x_{g,t-1}x_{g,t},\quad
&&y_{g,t}=p_{g,t-1}x_{g,t},\label{mono1}\\
&z_{g,t}=x_{g,t-1}p_{g,t},\quad
&&P_{g,t}=p^2_{g,t},\label{mono2}
\end{align}
\end{subequations}
hold true at optimality. 
To capture the binary requirement for commitment decisions, the following convex inequalities, that are referred to as ``McCormick constraints'', are enforced:
\begin{align}	
\max\{0, x_{g,t-1}+x_{g,t}-1\} \leq u_{g,t} \leq \min\{x_{g,t-1},x_{g,t}\}.\label{relax_mcc}
\end{align}
Now, constraint [\ref{cons_cost}] can be cast in the following linear form, with respect to the auxiliary variables:
\begin{align}	
	&c_{g,t} = \alpha_g o_{g,t} + \beta_g p_{g,t} +\nonumber\\ 
	&\quad\quad\;\;\,\gamma_g x_{g,t} + \gamma_g^{\uparrow} (x_{g,t}-u_{g,t}) + \gamma_g^{\downarrow}	 (x_{g,t-1}-u_{g,t}).\label{relax_cost}
\end{align}
Finally, we relax the nonconvex equations [\ref{mono1}]--[\ref{mono2}] with the following conic constraints: \\
\begin{subequations}\label{conic0}
	\noindent\begin{minipage}{0.422\hsize}
		\begin{equation}\label{conic1}
		\!\!\begin{bmatrix}
		x_{g,t}\! &\!\! u_{g,t}\! &\!\! y_{g,t}\!   \\
		u_{g,t}\! &\!\! x_{g,t-1}\! &\!\! p_{g,t-1}\!  \\
		y_{g,t}\! &\!\! p_{g,t-1}\! &\!\! o_{g,t-1}\!  
		\end{bmatrix}\!\succeq 0,
		\notag
		\addtocounter{equation}{1}
		\end{equation}
	\end{minipage}
	\begin{minipage}{0.262\hsize}
		\begin{equation} \label{conic2}            
		\notag
		\begin{bmatrix}
		x_{g,t-1}\! &\!\! u_{g,t}\! &\!\! z_{g,t}   \\
		u_{g,t}\! &\!\! x_{g,t}\! &\!\! p_{g,t}  \\
		z_{g,t}\! &\!\! p_{g,t}\! &\!\! o_{g,t}  
		\end{bmatrix}\!\succeq 0,
		\end{equation}
	\end{minipage}\begin{minipage}{0.3\hsize}
		\vspace{2.9mm}
		\hfill\;\;\normalfont[\begin{NoHyper}\ref{conic1},\ref{conic2}\end{NoHyper}]
	\end{minipage}
\end{subequations}\vspace{1mm}

\noindent
as well as a number of linear inequalities that are stated next.
\begin{definition}
	For each $g\in\mathcal{G}$, define $\mathcal{U}^{\mathrm{TSDP}}_g\subset\mathbb{R}^{T\times 4}$ to be the set of all quadruplets
	$(\mathbf{x}^{\top}_{g,\bullet},\mathbf{p}^{\top}_{g,\bullet},\mathbf{q}^{\top}_{g,\bullet},\mathbf{c}^{\top}_{g,\bullet})$, for which there exists $(\mathbf{u}^{\top}_{g,\bullet}, \mathbf{y}^{\top}_{g,\bullet}, \mathbf{z}^{\top}_{g,\bullet}, \mathbf{o}^{\top}_{g,\bullet})\in\mathbb{R}^{T\times4}$, such that for every $t\in\mathcal{T}$, the following constraints hold true:\vspace{0mm}
	\begin{itemize}
		\item[i)] 	The linear inequalities {\rm[\ref{cons_cap}]}, {\rm[\ref{cons_min}]}, {\rm[\ref{cons_ramp}]},
		\item[ii)]  The conic and linear constraints {\rm[\ref{relax_mcc}]}, {\rm[\ref{relax_cost}]}, {\rm[\ref{conic1}]}, {\rm[\ref{conic2}]} and {\rm[\ref{cons_K}]},
		\item[iii)] The linear inequalities {\rm[\ref{cons_Kd}]}, if $m^{\downarrow}_g>1$,
		\item[iv)]  The linear inequalities {\rm[\ref{cons_Ku}]}, if $m^{\uparrow}_g>1$.
	\end{itemize}
\end{definition}
Notice that for each $g\in\mathcal{G}$, the relaxed feasible set $\mathcal{U}^{\mathrm{TSDP}}_g$ is defined, by means of conic and linear inequalities that are convex. The validity of these inequalities is proven in \href{url}{SI}~\href{url}{Text}.

The definition of  $\mathcal{U}^{\mathrm{TSDP}}_g$ involves $2\times T$  third-order semidefinite constraints that can be enforced efficiently. Additionally, the overall number of inequalities grows linearly with respect to $T$, which is an improvement upon existing methods. On the other hand, the ramp and minimum up \& down constraints are incorporated into the valid inequalities, the present convex relaxation offers more accurate bounds, in the case of severe load variations.

\subsection{Relaxation of Network Constraints}
A state of the art method, given in \cite{madani2016promises}, for convex relaxation of AC power flow equations incorporates an auxiliary matrix variable $\mathbf{W}_t\in\mathbb{H}_n$, for each $t\in\mathcal{T}$, accounting for $\mathbf{v}_{\bullet, t}\mathbf{v}^{\ast}_{\bullet, t}$. Using the matrix $\mathbf{W}_t\in\mathbb{H}_n$, the AC network constraints [\ref{cons_AC1}]--[\ref{cons_AC4}] can be convexified as follows:
\begin{subequations}\label{relax_AC}\begin{align}
		\mathbf{d}_{t}+\mathrm{diag}\{\mathbf{W}_t\;\mathbf{Y}^{\ast}_t\}&=\mathbf{C}^{\top}(\mathbf{p}_{\bullet t}+i\mathbf{q}_{\bullet t})\label{relax_AC1}\\
		\lvert\mathrm{diag}\{\vec{\mathbf{C}}_t\;\mathbf{W}_t\;\vec{\mathbf{Y}}_t^{\ast}\}\rvert&\leq \mathbf{f}_{\mathrm{max};t}\label{relax_AC2}\\
		\lvert\mathrm{diag}\{\cev{\mathbf{C}}_t\;\mathbf{W}_t\;\cev{\mathbf{Y}}_t^{\ast}\}\rvert&\leq \mathbf{f}_{\mathrm{max};t}\label{relax_AC3}\\
		\mathbf{v}_{\mathrm{min}}\leq\lvert\mathrm{diag}\{\mathbf{W}_t\}\rvert&\leq \mathbf{v}_{\mathrm{max}}\label{relax_AC4}
\end{align}\end{subequations}
In formulation [\ref{relax_AC1}]--[\ref{relax_AC4}], the structure of matrix $\mathbf{W}_t$, i.e., 
\begin{align}
\mathbf{W}_t=\mathbf{v}_{\bullet, t}\mathbf{v}^{\ast}_{\bullet, t}
\label{W_fac}
\end{align}
is ignored, to make the model polynomially solvable. To remedy the absence of the non-convex equation [\ref{W_fac}], the relaxation can be strengthened through a combination of conic constraints, with the aim of enforcing the relation between $\mathbf{W}_t$ and $\mathbf{v}_{\bullet, t}$, implicitly. 

Observe that an arbitrary matrix $\mathbf{W}_t\in\mathbb{H}_n$ can be factored to $\mathbf{v}_{\bullet, t}\mathbf{v}^{\ast}_{\bullet, t}$, if and only if, it is rank-one and positive semidefinite:
\begin{align}
\mathrm{rank}\{\mathbf{W}_t\} = 1
\quad \wedge \quad 
\mathbf{W}_t\succeq 0.\label{rank_psd}
\end{align}
Although a rank constraint on $\mathbf{W}_t$ cannot be enforced efficiently, employing the convex constraint $\mathbf{W}_t\succeq 0$ leads to the semidefinite programming (SDP) relaxation of AC network constraints. For larger-scale systems, a graph-theoretic analysis divides the set of buses into several overlapping subsets $\mathcal{A}_1,\mathcal{A}_2,\ldots,\mathcal{A}_A\subseteq\mathcal{V}$, such that the relaxation could be represented with smaller conic constraints:
\begin{align} \label{relax:sdp}
\mathbf{W}_t[\mathcal{A}_k,\mathcal{A}_k]\succeq 0,\qquad\forall k\in\{1,2,\ldots,A\},
\end{align}
where for every $\mathcal{X}\subseteq\{1,\ldots,n\}$, the notation $\mathbf{W}_t[\mathcal{X},\mathcal{X}]$ represents the $|\mathcal{X}|\times|\mathcal{X}|$ principal submatrix of $\mathbf{W}_t$, whose rows and columns are chosen from $\mathcal{X}$. Choosing 
$\mathcal{A}_1,\mathcal{A}_2,\ldots,\mathcal{A}_A$,
based on the bags of an arbitrary tree decomposition of the network, leads to an equivalent but more efficient SDP relaxation \cite{madani2016promises}. A weaker, but far more tractable approach is the second-order cone programming (SOCP) relaxation which uses conic constraints of the following form:
\begin{align} \label{relax:socp}
	\begin{bmatrix}
		W_{k_1,k_1}\! &\!\! W_{k_1,k_2}\!  \\
		W_{k_2,k_1}\! &\!\! W_{k_2,k_2}\! 
	\end{bmatrix}\!\succeq 0,\qquad \forall (k_1,k_2)\in\mathcal{E}_t
\end{align}

To achieve a better balance between the strength of the convex relaxation and scalability, in this paper,
we use a third-order semidefinite programming (TSDP) relaxation which is described as follows:
%
\begin{definition}
For each $t\in\mathcal{T}$, let $\mathcal{N}^{\mathrm{TSDP}}_t\subset\mathbb{R}^{G\times 2}$ be the set of pairs $(\mathbf{p}^{\top}_{g,\bullet},\mathbf{q}^{\top}_{g,\bullet})$, for which there exists a Hermitian matrix $\mathbf{W}_t\in\mathbb{H}_n$, that satisfies constraints {\rm[\ref{relax_AC1}]}--{\rm[\ref{relax_AC4}]} and the following third-order semidefinite constraints
\begin{align}
\begin{bmatrix}
W_{k_1,k_1}\! &\!\! W_{k_1,k_2}\! &\!\! W_{k_1,k_3}   \\
W_{k_2,k_1}\! &\!\! W_{k_2,k_2}\! &\!\! W_{k_2,k_3}  \\
W_{k_3,k_1}\! &\!\! W_{k_3,k_2}\! &\!\! W_{k_3,k_3}  
\end{bmatrix}\!\succeq 0,
\end{align}
for every $(k_1,k_2,k_3)\in\bigcup_{k=1}^{A} \mathcal{A}_k\times\mathcal{A}_k\times\mathcal{A}_k$,
where 
$\mathcal{A}_1,\mathcal{A}_2,\ldots,\mathcal{A}_A\subseteq\mathcal{V}$,
are the bags associated with an arbitrary tree decomposition of the network $\mathcal{H}_t$, for which $|\mathcal{A}_k|\geq 3$, $\forall k \in\{1,\ldots A\}$.\label{def3}
\end{definition}

\begin{figure}
	\centering
	\includegraphics[width=0.82\linewidth]{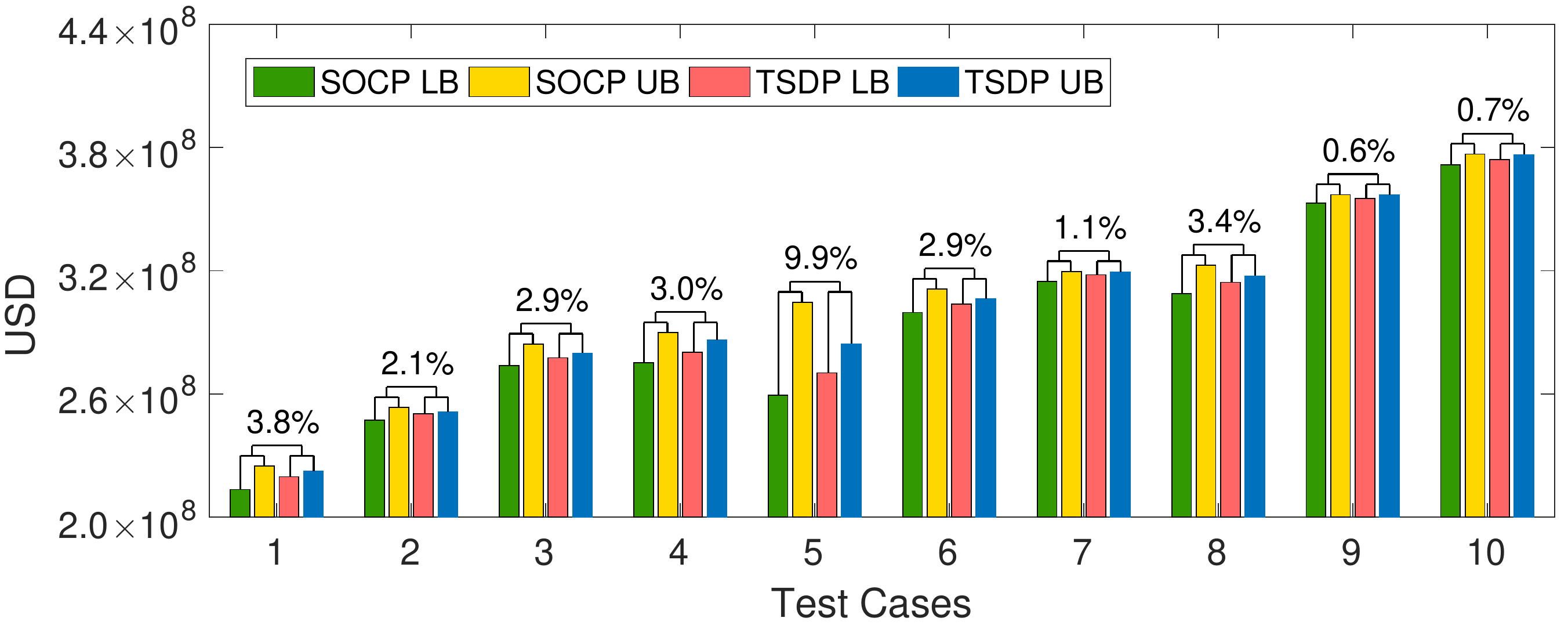}
	\caption{Comparisons between the performance of TSDP and SOCP relaxations of the AC network constraints on ten instances of PEGASE 1354-bus system. The difference between optimality gaps are shown.}\label{table_res_OPF}\vspace{-2.0mm}
\end{figure}

\section*{Comparison of the models}
We have compared the performance of the proposed TSDP relaxation 
(Definition \ref{def3}) with
the SOCP relaxation given by [\ref{relax:socp}] 
in terms of: i) upper bound (UB) given by the cost of the feasible solution returned and; ii) convex lower bound (LB) on the optimal cost. To that end, ten instances based on each European system 
in Table \ref{table_res_AC}  are considered. 
In all instances of the largest PEGASE 2869-bus, 9241-bus and 13659-bus cases, the Newton-Raphson's local search method successfully converges to a feasible operating point when started with an initial point from TSDP relaxation (discussed in detail in Section Recovering a Feasible Solution); however, it fails to converge when started with points from the SOCP relaxation. For the instances on the smallest PEGASE 1354-bus system, both SOCP and TSDP relaxations lead to feasible solutions through Newton-Raphson's local search. 
Figure \ref{table_res_OPF} displays the optimal relaxation objective value (LB) and the cost of the resulting feasible solution (UB) for this system. Although very fast to solve, the SOCP relaxation is substantially worse than the TSDP relaxation, as seen in Figure \ref{table_res_OPF}.
Motivated by the aforementioned observation, we propose the TSDP relaxation to convexify the AC network equations [\ref{cons_AC1}]--[\ref{cons_AC4}], in power system scheduling.

\section*{Discussion and Conclusions}
In this paper, we study the problem of optimizing grids operation throughout a planning horizon, based on the available resources for supply and transmission of electricity. This fundamental problem is heavily investigated for decades and need to be solved on a daily basis by independent system operators and utility companies. The challenge is twofold: first, determining a massive number of highly correlated binary decisions that account for commitment of generators; secondly, finding the most economic transmission strategy in accordance with laws of physics and technological limitations. 

We propose a third-order semidefinite programming (TSDP) method that is equipped with an accurate physical model for the flow of electricity and offers massive scalability in the number of generating units and grid size. While co-optimization of supply and transmission under the full physical model has long been put forward as a direction to boost the efficiency and reliability of operation, scalability has been the main bottleneck to-date. Significant improvement over the state-of-the-art methods is validated on day-ahead grid scheduling problems, on the largest publicly available real-world data.

Given the simplicity of the linear algebraic operations on 3x3 Hermitian matrices, a direction of interest is to build a highly parallel numerical algorithm for solving large-scale TSDP problems on graphical processing units high-performance computing facilitates.

\section*{Materials and Methods}	
	This section details the procedure for recovering a feasible solution to the scheduling from the optimal solution of the relaxed problem [\ref{PSS_relax_obj}]--[\ref{cons_network_relax}] and data generation.
	
\subsection*{Recovering a Feasible Solution}

Let $(
\mathbf{x}^{\mathrm{opt}},
\mathbf{p}^{\mathrm{opt}},
\mathbf{q}^{\mathrm{opt}},
\mathbf{c}^{\mathrm{opt}}
)$ be an optimal solution to the relaxed problem [\ref{PSS_relax_obj}]--[\ref{cons_network_relax}]. 
If all entries of $\mathbf{x}^{\mathrm{opt}}$ turn out to be integer, and there exists a matrix 
$\mathbf{v}\in\mathbb{C}^{\mathcal{V}\times \mathcal{T}}$ that satisfies constraints [\ref{cons_AC1}]--[\ref{cons_AC4}], then $(
\mathbf{x}^{\mathrm{opt}},
\mathbf{p}^{\mathrm{opt}},
\mathbf{q}^{\mathrm{opt}},
\mathbf{c}^{\mathrm{opt}}
)$ is a globally optimal solution to problem [\ref{PSS_obj}]--[\ref{cons_network}].
However, the relaxation is often inexact, and solutions to the relaxed problem [\ref{PSS_relax_obj}]--[\ref{cons_network_relax}] are not necessarily feasible for problem [\ref{PSS_obj}]--[\ref{cons_network}]. In such cases, a recovery process is needed to transform $(
\mathbf{x}^{\mathrm{opt}},
\mathbf{p}^{\mathrm{opt}},
\mathbf{q}^{\mathrm{opt}},
\mathbf{c}^{\mathrm{opt}}
)$ to a feasible and near-optimal solution for problem [\ref{PSS_obj}]--[\ref{cons_network}].

\begin{algorithm}[t]
{\small
	\caption{\small Recovering a Feasible Solution}
	\label{alg}
	\begin{algorithmic}
		\Require{The optimal unit commitment solution $\mathbf{x}^{\mathrm{opt}}\in\mathbb{R}^{G\times T}$ to problem [\ref{PSS_relax_obj}]--[\ref{cons_network_relax}]:}	\vspace{1mm}
		\For{$g=1\ldots,G$}\vspace{1mm}
		\For{$t=1-\max\{m^{\uparrow}_g,m^{\downarrow}_g\},\ldots,0$}\vspace{1mm}
		\State Set $x^{\mathrm{feas}}_{g,t}$ according to the initial state of unit $g$.
		\EndFor
		\EndFor
		\For{$g=1\ldots,G$}\vspace{1mm}
		\For{$t=1\ldots,T$}\vspace{1mm}
		\State $a^{\uparrow}\!\!\gets\!\max\{x^{\mathrm{feas}}_{g,\tau}-x^{\mathrm{feas}}_{g,\tau-1}\;|\;\forall\tau\! \in\!\{t-m^{\uparrow}_g+1,\ldots,t-1\}\}$\vspace{1mm}
		\State $a^{\downarrow}\!\!\gets\!\max\{x^{\mathrm{feas}}_{g,\tau-1}-x^{\mathrm{feas}}_{g,\tau}\;|\;\forall\tau\! \in\!\{t-m^{\downarrow}_g+1,\ldots,t-1\}\}$\vspace{1mm}
		\If{$a^{\uparrow}=a^{\downarrow}=1$}
		\State Declare failure.
		\Else
		\If{$a^{\uparrow}=1$}
		\State $x^{\mathrm{feas}}_{g,t}\gets 1$
		\EndIf
		\If{$a^{\downarrow}=1$}
		\State $x^{\mathrm{feas}}_{g,t}\gets 0$
		\EndIf
		\If{$a^{\uparrow}=a^{\downarrow}=0$}
		\State $x^{\mathrm{feas}}_{g,t} \gets\mathrm{round}\{x^{\mathrm{opt}}_{g,t}+0.25\} $
		\EndIf
		\EndIf
		\EndFor
		\EndFor
		\State \Return $\mathbf{x}^{\mathrm{feas}}$
	\end{algorithmic}}
\end{algorithm}

As demonstrated by Table \ref{table_res_AC}, on average, only a small portion of the binary variables remain fractional after solving the proposed TSDP relaxation problem. 
In all of our experiments, a feasible candidate for $\mathbf{x}$ is obtained, through the Algorithm \ref{alg}, which simply rounds each entry of $\mathbf{x}^{\mathrm{opt}}$ subject to minimum up and down time constraints [\ref{cons_min_l}] and [\ref{cons_min_u}]. 

Another challenge is finding a feasible voltage profile 
$
\mathbf{v}=[
\mathbf{v}_{\bullet, 1}|
\mathbf{v}_{\bullet, 2}|
\ldots|
\mathbf{v}_{\bullet, T}
]\in\mathbb{C}^{\mathcal{V}\times \mathcal{T}},
$
based on a solution to the relaxed problem [\ref{PSS_relax_obj}]--[\ref{cons_network_relax}]. 
If the rank constraints [\ref{rank_psd}] are not satisfied at optimality, then the relaxation of network equations is not exact and it is not possible to factorize the resulting matrices 
$\mathbf{W}^{\mathrm{opt}}_1,
\mathbf{W}^{\mathrm{opt}}_2,\ldots,$ $
\mathbf{W}^{\mathrm{opt}}_T$, 
in the form of equation [\ref{W_fac}]. A ``recovery algorithm'' is introduced in \cite{madani2016promises}, for finding an approximate vector of voltages $\hat{\mathbf{v}}_{\bullet, t}$ based on $\mathbf{W}^{\mathrm{opt}}_t$, which minimizes the overall mismatch (i.e., violation of network equations). In order to obtain voltage profiles with no mismatch, we feed the outcome of the recovery algorithm from \cite{madani2016promises} as the initial point to Newton-Raphson's local search algorithm. This procedure is described next:
\begin{enumerate}
	\item Find a feasible matrix of commitment decisions $\mathbf{x}^{\mathrm{feas}}$ via Algorithm~\ref{alg}. 
	\item For every $t=1\ldots,T$:
	\begin{enumerate}
		\item[i)] Obtain an approximate voltage profile $\hat{\mathbf{v}}_{\bullet, t}$ from $\mathbf{W}^{\mathrm{opt}}_t$ based on the recovery algorithm in \cite{madani2016promises}. 
		\item[ii)] Start with
		$\mathbf{p}^{\mathrm{opt}}_{\bullet, t}$,
		$\mathbf{q}^{\mathrm{opt}}_{\bullet, t}$ and
		$\hat{\mathbf{v}}_{\bullet, t}$, as the initial point to
		search locally for a triplet of vectors
		$\mathbf{p}^{\mathrm{feas}}_{\bullet, t}\in\mathbb{R}^G$,
		$\mathbf{q}^{\mathrm{feas}}_{\bullet, t}\in\mathbb{R}^G$ and
		$\mathbf{v}^{\mathrm{feas}}_{\bullet, t}\in\mathbb{C}^{\mathcal{V}}$, 
		that minimizes the objective function
		$\sum^{G}_{g=1}{\alpha_{g} p^2_{g,t} + \beta_{g} p_{g,t}}$, subject to the constraints [\ref{cons_AC1}]--[\ref{cons_AC4}] and
		\begin{subequations}
		\begin{align}
	&\ubar{q}_g x^{\mathrm{feas}}_{g,t} \leq  q_{g,t} \leq \bar{q}_g x^{\mathrm{feas}}_{g,t}\\ 
	&\ubar{p}_g x^{\mathrm{feas}}_{g,t} \leq  p_{g,t} \leq \bar{p}_g x^{\mathrm{feas}}_{g,t}\\ 
	&p_{g,t}\geq p^{\mathrm{feas}}_{g,t-1} - r_g x^{\mathrm{feas}}_{g,t} - s_g (1-x^{\mathrm{feas}}_{g,t}) \\
	&p_{g,t}\leq p^{\mathrm{feas}}_{g,t-1} + r_g x^{\mathrm{feas}}_{g,t-1} + s_g (1-x^{\mathrm{feas}}_{g,t-1}).
		\end{align}
		\end{subequations}
		\item[iii)]
		Derive the feasible cost values $c^{\mathrm{feas}}_{g,t}$, according to the equation [\ref{cons_cost}].
	\end{enumerate}
	\item Report $(\mathbf{x}^{\mathrm{feas}},\mathbf{p}^{\mathrm{feas}},\mathbf{q}^{\mathrm{feas}})$ as the output schedule/dispatch and $\mathbf{v}$ as the corresponding voltage profile. The following quantity serves as an upperbound for relative distance from global optimality:
	\begin{align}
	\!\!\mathrm{Gap}\leq 100\times\frac{
		\sum_{t=1}^{T}{\sum_{g=1}^{G}{(c^{\mathrm{feas}}_{g,t}-c^{\mathrm{opt}}_{g,t})}}}
	{\sum_{t=1}^{T}{\sum_{g=1}^{G}{c^{\mathrm{feas}}_{g,t}}}}
	\end{align}
\end{enumerate}
We have used the procedure described above for the experiments presented in Table \ref{table_res_AC}, and in all cases, a feasible solution could be found within the violation tolerance of the constraint ($10^{-5}$ per-unit).

\subsection*{Data Generation}
	
The network data for IEEE and European systems is obtained from the MATPOWER package \cite{zimmerman2011matpower,josz2016ac}. Hourly load changes for the day-ahead at all buses are considered proportional to the numbers reported in \cite{khodaei2010transmission}.
In each experiment, the cost coefficients 
$\alpha_g$, $\beta_g$, $\gamma_g$, $\gamma^{\downarrow}_g$ and $\gamma^{\uparrow}_g$
are chosen uniformly between zero and $~1~\$/(\mathrm{MW.h})^2$, $~10~\$/(\mathrm{MW.h})$, $~100~\$$, $~30~\$$ and $~50~\$$, respectively. 
The ramp limits of each generating unit are set to 
$r_g = s_g = \max\{\bar{p}_g/4,\ubar{p}_g\}$. For each generating unit, the minimum up and down limits $m^{\uparrow}_g$ and $m^{\downarrow}_g$ are randomly selected in such a way that $m^{\uparrow}_g-1$ and $m^{\downarrow}_g-1$ have Poisson distribution with parameter $4$. The initial status of generators at time period $t=0$ is found by solving a single period economic dispatch problem corresponding to the demand at time $t=1$. For each generating unit $g\in\mathcal{G}$, it is assumed that the initial status has been maintained exactly since time period $t=-t^{(0)}_g$, where $t^{(0)}_g$ has Poisson distribution with parameter $4$. For simplicity, all of the generating units with negative capacity are removed.
All simulations are run in MATLAB using a workstation with an Intel 3.0 GHz, 12-core CPU, and 256 GB RAM. The CVX package version 3.0 \cite{cvx} and MOSEK version 8.0 \cite{mosek} are used for solving semidefinite programming problems. The data set as well as the log files of the optimization runs are available for download at:
\href{http://ieor.berkeley.edu/~atamturk/data/tsdp}{http://ieor.berkeley.edu/$\sim$atamturk/data/tsdp}.


\subsection*{Acknowledegement}
A. Atamt\"urk is supported, in part, 
	by grant FA9550-10-1-0168 from the Office of the Assistant Secretary of Defense for Research and Engineering. 



\bibliographystyle{IEEEtran}
\bibliography{PSS}

\pagebreak

\section*{Appendix} 
%
%
%
%
%


A proof of validity for conic and linear inequalities [\ref{conic0}] and [\ref{cons_KK}] is provided in this section.



\begin{proposition}
	Inequalities {\rm[\ref{conic1}]}, {\rm[\ref{conic2}]} and {\rm[\ref{cons_K}]} are valid for every pair of vectors $(\mathbf{x},\mathbf{p})\in\mathbb{R}^{G\times T}$ that satisfy constraints
	{\rm[\ref{cons_bin}]}, {\rm[\ref{cons_cap_p}]}, {\rm[\ref{cons_ramp_pu}]}, {\rm[\ref{cons_ramp_pl}]}, {\rm[\ref{cons_min_u}]} and {\rm[\ref{cons_min_l}]}. Additionally, if 
	$m^{\downarrow}_g \geq 2$, then constraint {\rm[\ref{cons_Kd}]} and if $m^{\uparrow}_g \geq 2$, then constraint {\rm[\ref{cons_Ku}]} is valid.
\end{proposition}
\textit{\textbf{Proof:}}
For every 
$(g,t)\in\{1,\ldots,G\}\times\{1,\ldots,T\}$, define the vector of monomials:
\begin{align}
\boldsymbol{\delta}_{g,t} \triangleq [&x_{g,t-1},\; x_{g,t},\; p_{g,t-1},\; p_{g,t},\nonumber\\ 
&\ubar{h}_{g,t-1},\; \bar{h}_{g,t-1},\; \ubar{h}_{g,t},\; \bar{h}_{g,t},\; \ubar{a}_{g,t},\; \bar{a}_{g,t},\nonumber\\
&x_{g,t-1}\ubar{h}_{g,t-1},\; x_{g,t-1}\bar{h}_{g,t-1},\; x_{g,t-1}\ubar{h}_{g,t},\nonumber\\ &x_{g,t-1}\bar{h}_{g,t},x_{g,t-1}\ubar{a}_{g,t},\; x_{g,t-1}\bar{a}_{g,t},\nonumber\\
&x_{g,t}\ubar{h}_{g,t-1},\; x_{g,t}\bar{h}_{g,t-1},\; x_{g,t}\ubar{h}_{g,t},\nonumber\\ 
&x_{g,t}\bar{h}_{g,t},\; x_{g,t}\ubar{a}_{g,t},\; x_{g,t}\bar{a}_{g,t}]^{\top},
\end{align}
where
\begin{subequations}
	\begin{align}
	\ubar{h}_{g,t} &\triangleq  \sqrt{p_{g,t}- \ubar{p}_g x_{g,t}}\;,\qquad 	
	\bar{h}_{g,t} \triangleq \sqrt{\bar{p}_g x_{g,t} - p_{g,t}}\;,  \\
	\ubar{a}_{g,t} &\triangleq \sqrt{s_g + (r_g - s_g) x_{g,t-1} + p_{g,t-1} - p_{g,t} }\;, \\ 
	\bar{a}_{g,t} &\triangleq \sqrt{s_g + (r_g - s_g) x_{g,t} - p_{g,t-1} + p_{g,t} }\;.
	\end{align}
\end{subequations}
%
Define $\boldsymbol{\Delta}_{g,t}$ as the $22\times22$ symmetric matrix formed by multiplying $\boldsymbol{\delta}_{g,t}$ by its transpose:
\begin{align}
\boldsymbol{\Delta}_{g,t}\triangleq\boldsymbol{\delta}_{g,t}\boldsymbol{\delta}_{g,t}^{\top}.
\end{align}
Observe that $\boldsymbol{\Delta}_{g,t}$ is positive semidefinite, and as a consequence, every principle submatrix of $\boldsymbol{\Delta}_{g,t}$ is positive semidefinite, as well. Considering submatrices $\boldsymbol{\Delta}_{g,t}[\{2,1,3\},\{2,1,3\}]$ and
$\boldsymbol{\Delta}_{g,t}[\{1,2,4\},\{1,2,4\}]$, concludes the conic constraints [\ref{conic1}] and [\ref{conic2}], respectively. Moreover, the constraint [\ref{cons_K}] encapsulates $24$ linear inequalities, and it is straightforward to verify that inequalities $k$, $k+6$, $k+12$ and $k+18$ can be concluded from
\begin{align}
\!\!\boldsymbol{\Delta}_{g,t}[\{k+4,k+10,k+16\},\{k+4,k+10,k+16\}]\succeq 0,\!
\end{align}
for each $k=1,2,\ldots,6$. This completes the proof of [\ref{cons_K}].

In order to prove the validity constraint [\ref{cons_Kd}], suppose that $m^{\downarrow}_g \geq 2$, and consider the following vector of monomials:
\begin{align}
\boldsymbol{\phi}^{\downarrow}_{g,t} \triangleq [&w^{\downarrow}_{g,t},\;
x_{g,t-1}w^{\downarrow}_{g,t},\;
\ubar{h}_{g,t-1}w^{\downarrow}_{g,t},\;
\bar{h}_{g,t-1}w^{\downarrow}_{g,t}]^{\top},
\end{align}
where
\begin{align}
w^{\downarrow}_{g,t} \triangleq \sqrt{1-x_{g,t-2}+x_{g,t-1}-x_{g,t}}\;.
\end{align}
Observe that all four inequalities encapsulated in [\ref{cons_Kd}] can be concluded from the conic inequality $\boldsymbol{\phi}^{\downarrow}_{g,t}(\boldsymbol{\phi}^{\downarrow}_{g,t})^{\top}\succeq0$.

If $m^{\uparrow}_g \geq 2$, the validity of the constraint [\ref{cons_Ku}] can be similarly proven by defining
\begin{align}
w^{\uparrow}_{g,t} \triangleq\sqrt{x_{g,t-2}-x_{g,t-1}+x_{g,t}}\;,
\end{align}
and forming the vector of monomials
\begin{align}
\boldsymbol{\phi}^{\uparrow}_{g,t} \triangleq [&w^{\uparrow}_{g,t},\;
x_{g,t-1}w^{\uparrow}_{g,t},\;
\ubar{h}_{g,t-1}w^{\uparrow}_{g,t},\;
\bar{h}_{g,t-1}w^{\uparrow}_{g,t}]^{\top}.
\end{align}
Finally, the four inequalities from [\ref{cons_Ku}] can be inferred from the conic inequality $\boldsymbol{\phi}^{\uparrow}_{g,t}(\boldsymbol{\phi}^{\uparrow}_{g,t})^{\top}\succeq0$.
\qed

\end{document}